\documentclass[12pt,reqno]{article}

\usepackage[usenames]{color}
\usepackage{amssymb}
\usepackage{amsmath}
\usepackage{amsthm}
\usepackage{amsfonts}
\usepackage{amscd}
\usepackage{graphicx}

\usepackage[colorlinks=true,
linkcolor=webgreen,
filecolor=webbrown,
citecolor=webgreen]{hyperref}

\definecolor{webgreen}{rgb}{0,.5,0}
\definecolor{webbrown}{rgb}{.6,0,0}

\usepackage{color}
\usepackage{fullpage}
\usepackage{float}

\newcommand{\seqnum}[1]{\href{https://oeis.org/#1}{\rm \underline{#1}}}

\begin{document}

\def\refname {References}
\newtheorem{theorem} {Theorem}[section]
\newtheorem{lemma}[theorem] {Lemma}
\newtheorem{definition}[theorem] {Definition}
\newtheorem{example}[theorem] {Example}
\newtheorem{corollary}[theorem] {Corollary}
\newtheorem{proposition} [theorem]{Proposition}
\newtheorem{problem}{Problem}

\def\bibname {\bf References}

\title{\bf On Sum of a Polynomial Multiplied by Generalized Fibonacci Numbers}
\date{}
\author{Ivan Hadinata (updated 17 oktober 2024)\\\\
Department of Mathematics,\\ Faculty of Mathematics and Natural Sciences, 
Gadjah Mada University,\\ Yogyakarta, Indonesia\\
\href{mailto: ivanhadinata2005@mail.ugm.ac.id}{\tt ivanhadinata2005@mail.ugm.ac.id}}
\begingroup
\let\newpage\relax
\maketitle
\endgroup
\begin{abstract}
     Given that $a,b\in\mathbb N$, $c_0,c_1\in\mathbb Z$, $(c_0,c_1)\neq (0,0)$, and a generalized Fibonacci sequence $(s_n)_{n\geq 0}$ where $s_0 = c_0$, $s_1 = c_1$, and $s_{n+1}=as_{n}+bs_{n-1}$ for all positive integers $n$. In this paper, we get the result that for every polynomials $P(x)$ with real coefficients, we can always find three polynomials $F_1(x), G_1(x), H_1(x)$ (not necessarily distinct) with real coefficients satisfying the identity: $\;2\sum_{k=1}^{n}P(k)s_{k-1} = F_1(n)s_{n+1} + G_1(n)s_n + H_1(n), \;\forall n\in\mathbb N$. Furthermore, we serve two constraints for $(s_n)_{n\geq 0}$: one constraint implies that there are infinitely many triples $(F_1(x), G_1(x), H_1(x))$ satisfying the identity $\;2\sum_{k=1}^{n}P(k)s_{k-1} = F_1(n)s_{n+1} + G_1(n)s_n + H_1(n), \;\forall n\in\mathbb N$, while another constraint implies that there is only one triple $(F_1(x), G_1(x), H_1(x))$ satisfying the identity $\;2\sum_{k=1}^{n}P(k)s_{k-1} = F_1(n)s_{n+1} + G_1(n)s_n + H_1(n), \;\forall n\in\mathbb N$.
\end{abstract}

\section{Introduction} 
Fibonacci sequence \cite{koshy with fibonacci and lucas}, Lucas sequence \cite{koshy with fibonacci and lucas}, Pell sequence \cite{koshy with pell}, and Jacobsthal sequence \cite{sloane} are well-known examples of second order linear recurrence sequences. All of them can be generalized by sequence $(s_n)_{n\geq 0}$ with $\, s_0=c_0$, $\, s_1=c_1$, $\, s_{n+1} = as_{n} + bs_{n-1}, \;\forall n\in\mathbb N \:$ where $\: a,b\in\mathbb N;\: c_0,c_1\in\mathbb Z$, $(c_0,c_1)\neq (0,0)$. Indeed, the sequence $(s_n)_{n\geq 0}$ is a special case of sequence $\{w_n\}$ defined by Horadam in \cite{horadam}. In this case, $(s_n)_{n\geq 0}$ can be rewritten as Horadam sequence $\,\{w_n(c_0, c_1; a, -b)\}\,$ where $\: a,b\in\mathbb N;\: c_0,c_1\in\mathbb Z$, $(c_0,c_1)\neq (0,0)$. If $(a,b,c_0,c_1)=(1,1,0,1)$, then $(s_n)_{n\geq 0}$ becomes the Fibonacci sequence $(F_n)_{n\geq 0}$. If $(a,b,c_0,c_1)=(1,1,2,1)$, then $(s_n)_{n\geq 0}$ becomes the Lucas sequence $(L_n)_{n\geq 0}$. If $(a,b,c_0,c_1)=(2,1,0,1)$, then $(s_n)_{n\geq 0}$ becomes the Pell sequence $(P_n)_{n\geq 0}$. If $(a,b,c_0,c_1) = (1,2,1,0)$, then $(s_n)_{n\geq 0}$ becomes the Jacobsthal sequence $(J_n)_{n\geq 0}$.
\begin{definition}
    Let $(s_n)_{n\geq 0}$ be a generalized Fibonacci sequence defined by:
    \begin{equation*}
        s_0 = c_0, \;\; s_1 = c_1, \;\; \textrm{and} \;\;\;(\forall n\in\mathbb N)\; s_{n+1} = as_n + bs_{n-1}
    \end{equation*}
    with $\; a,b\in\mathbb N$; $\;\; c_0,c_1\in\mathbb Z$; $\; (c_0,c_1)\neq (0,0)$.
\end{definition}
\vspace{0.5cc} 

In this paper, we would like to study about summation $\: \sum_{k=1}^nP(k)s_{k-1} \:$ where $P(x)$ is a polynomial with real coefficients. This form is a development of the series like $\: \sum_{k=1}^{n}k^dF_{k-1}$, $\: \sum_{k=1}^{n}k^dL_{k-1}$, $\: \sum_{k=1}^{n}k^dP_{k-1}$, $\: \sum_{k=1}^{n}k^dJ_{k-1}$, and many others (note: $d\in\mathbb N_0$). Some articles have discussed the related kinds of these summations, for example \cite{g. ledin} and \cite{bro alfred}. Ledin \cite{g. ledin} proposes an identity for summation $\sum_{k=1}^{n}k^m F_k$ with $m\in\mathbb N_0$. He presents that $\sum_{k=1}^{n}k^m F_k = P_2(m,n)F_{n+1} + P_1(m,n)F_n + C(m)$ where $P_1(m,n)$ and $P_2(m,n)$ are polynomials in $m$ of degree $n$, and $C(m)$ is a real constant depending on $m$. Besides that, in \cite{bro alfred}, Alfred discovers the general formula of summation $\sum_{k=1}^{n}k^mF_{k+r}$ with $m,r\in\mathbb N$ by using a finite difference approach.

\vspace{0.3cc}

We are also firstly motivated by problem $1410$ Spring $2024$ of Pi Mu Epsilon journal \cite{pi mu epsilon} which is proposed by Kenny B. Davenport. The problem considers Pell sequence $(P_n)_{n\geq 0}$ where $P_0=0$, $P_1=1$, $P_{k+2} = 2P_{k+1}+P_k, \:\forall k\geq 0$. There are two identities that must be proved: 
\begin{equation}\label{eq1}
    2\sum_{k=1}^{n}kP_{k-1} = nP_{n+1} - (n+1)P_n; \quad 2\sum_{k=1}^{n}k^2P_{k-1} = (n^2+1)P_{n+1} - (n^2+2n)P_n -1.
\end{equation}
These identities can be proved easily by induction on $n$, but the more challenging question in this problem is "can we conjecture what shape of the identity obtained by the summation $2\sum_{k=1}^{n}k^dP_{k-1}$ ?". Luckily, we find an interesting identity of the summation $\: 2\sum_{k=1}^{n}k^dP_{k-1} \:$ for general $d\in\mathbb N_0$ and it generalizes two identities in \eqref{eq1}. This identity is in the following theorem.
\vspace{0.5cc}
\begin{theorem}\label{first theorem}
    Let $d\in\mathbb N_0$, $(P_n)_{n\geq 0}$ be Pell sequence, and $A_d = (a_{i,j})_{i,j=1}^{d+1}$ be a $(d+1)\times(d+1)$ square matrix defined by:
    \begin{equation*}
       a_{i,j} = \begin{cases}
            0 & \textrm{if} \quad i>j \\
            2 & \textrm{if} \quad i=j \\
            (2^{j-i}+2)\binom{j-1}{i-1} & \textrm{if} \quad i<j.
        \end{cases}
    \end{equation*}
    The matrix $A_d$ is invertible because its determinant is $2^{d+1}$ which is non-zero. Therefore there exists an ordered $(d+1)$-tuple $(\lambda_0(d), \lambda_1(d), \ldots, \lambda_d(d))\in\mathbb R^{d+1}$ which is a unique solution of matrix equation
\begin{equation*}
    A_d(\lambda_0(d), \lambda_1(d), \ldots, \lambda_d(d))^T = 2\left(\binom{d}{0}2^d, \binom{d}{1}2^{d-1}, \ldots, \binom{d}{d}2^0\right)^T.
\end{equation*}
    Then the polynomials $F(x)$, $G(x)$, $H(x)$ defined by:
    \begin{align*}
        F(x) = \sum_{i=0}^{d}\lambda_i(d)x^i, \; G(x) = -2(x+1)^d + \left(\sum_{i=0}^{d}(x+1)^i\lambda_i(d)\right), \; H(x) = 2^{d+1} - \sum_{i=0}^{d}(2^i+2)\lambda_i(d)
    \end{align*}
    satisfy the identity
    \begin{equation}\label{main identity in first theorem}
        (\forall n\in\mathbb N) \;\; 2\sum_{k=1}^{n}k^dP_{k-1} = F(n)P_{n+1} + G(n)P_n + H(n).
    \end{equation}
\end{theorem}
\vspace{0.5cc}

\begin{example}
    Some identities of $\: 2\sum_{k=1}^{n}k^dP_{k-1} \:$ for lower degrees $d$ are:
    \begin{align*}
        d=0 \;\implies\; 2\sum_{k=1}^{n}P_{k-1} =\, &P_{n+1} - P_n - 1,  \\
        d=1 \;\implies\; 2\sum_{k=1}^{n}kP_{k-1} =\, &nP_{n+1} - (n+1)P_n,  \\
        d=2 \;\implies\; 2\sum_{k=1}^{n}k^2P_{k-1} =\, &(n^2+1)P_{n+1} - (n^2+2n)P_n - 1,  \\
        d=3 \;\implies\; 2\sum_{k=1}^{n} k^3P_{k-1} =\, &(n^3 + 3n -3)P_{n+1} - (n^3+3n^2+1)P_n +3,  \\
        d=4 \;\implies\; 2\sum_{k=1}^{n}k^4P_{k-1} =\, &(n^4+6n^2-12n+13)P_{n+1} - (n^4+4n^3+4n-6)P_n - 13, \\
        d=5 \;\implies\; 2\sum_{k=1}^{n}k^5P_{k-1} =\, &(n^5 + 10n^3 -30n^2 +65n -75)P_{n+1} \\
        &- (n^5+5n^4+10n^2-30n +31)P_n + 75, \\
        d=6 \;\implies\; 2\sum_{k=1}^{n}k^6P_{k-1} =\, &(n^6 + 15n^4 -60n^3 +195n^2 -450n + 511)P_{n+1} \\
        &- (n^6 +6n^5 +20n^3 -90n^2+186n -210)P_n - 511.
    \end{align*}
\end{example}

\vspace{1cc}

Theorem \eqref{first theorem} can be proved by using induction and some other rigid computations. It can be generalized to be Theorem \eqref{matrix B and configuration of F,G,H} preceded by Lemma \eqref{lemma k^d s(k-1)}. As our attempt to generalize $(P_n)_{n\geq 0}$ to $(s_n)_{n\geq 0}$, we can generalize the identity \eqref{main identity in first theorem} to be identity \eqref{k^d s(k-1)}. Identity \eqref{k^d s(k-1)} states that $\: \quad 2\sum_{k=1}^{n}k^ds_{k-1} = F_{1,d}(n)s_{n+1} + G_{1,d}(n)s_n + H_{1,d}(n) \:$ for all positive integers $n$, where $F_{1,d}(x)$, $G_{1,d}(x)$, and $H_{1,d}(x)$ are certain polynomials with real coefficients. Later, an example of such polynomials $F_{1,d}(x)$, $G_{1,d}(x)$, $H_{1,d}(x)$ is stated in Theorem \eqref{matrix B and configuration of F,G,H}.

\vspace{0.3cc}

Finally, our goal is to present the main Theorems \eqref{main theorem 1}, \eqref{main theorem 2}, and \eqref{main theorem 3}. For every polynomials $P(x)\in\mathbb R[x]$, there exists polynomials with real coefficients $\,F_1(x)$, $\,G_1(x)$, $\,H_1(x)$ satisfying the identity $\:2\sum_{k=1}^{n}P(k)s_{k-1} = F_1(n)s_{n+1} + G_1(n)s_n + H_1(n)\:$ for all positive integers $n$. Later, we consider the sequence $(s_n)_{n\geq 0}$ into $2$ cases: (i) when both of $\: 2c_1 - (a+\sqrt{a^2+4b})c_0 \:$ and $\: 2c_1 - (a-\sqrt{a^2+4b})c_0 \:$ are non-zero, and (ii) when one of $\: 2c_1 - (a+\sqrt{a^2+4b})c_0 \:$ or $\: 2c_1 - (a-\sqrt{a^2+4b})c_0 \:$ equals zero. When both of $\; 2c_1 - (a+\sqrt{a^2+4b})c_0 \;$ and $\; 2c_1 - (a-\sqrt{a^2+4b})c_0 \;$ are non-zero, by the help of Theorem \eqref{polinom gamma1, gamma2, gamma3}, we ensure that the such triple $(F_1(x), G_1(x), H_1(x))$ is unique with dependent of $P(x)$. If we set $P(x) = a_mx^m + a_{m-1}x^{m-1}+\cdots + a_1x + a_0 \,$ (where $m\in\mathbb N_0$ and $a_0, a_1, \ldots, a_m\in\mathbb R$), we obtain $\: F_1(x) = \sum_{d=0}^{m}a_dF_{1,d}(x)$, $\: G_1(x) = \sum_{d=0}^{m}a_dG_{1,d}(x)$, $\: H_1(x) = \sum_{d=0}^{m}a_dH_{1,d}(x) \:$ where $F_{1,d}(x)$, $G_{1,d}(x)$, $H_{1,d}(x)$ are polynomials which stated in Theorem \eqref{matrix B and configuration of F,G,H}. On the other hand, when one of $\: 2c_1 - (a+\sqrt{a^2+4b})c_0 \:$ or $\: 2c_1 - (a-\sqrt{a^2+4b})c_0 \:$ equals zero, so for a single polynomial $P(x)\in\mathbb R[x]$, there are infinitely many triples $(F_1(x), G_1(x), H_1(x))$ in $\mathbb R[x]^3$ such that the identity $\:2\sum_{k=1}^{n}P(k)s_{k-1} = F_1(n)s_{n+1} + G_1(n)s_n + H_1(n), \:\forall n\in\mathbb N\:$ holds.

\vspace{1.5cc}

\section{Some Important Facts}\label{section 2}
\noindent
For more convenient, let us define $j_1$ and $j_2$ as follows.
\begin{definition}
    Define $\; j_1 = a+\sqrt{a^2+4b} \;$ and $\; j_2 = a-\sqrt{a^2+4b}$.
\end{definition}
\vspace{0.5cc}
\noindent
We consider that the sequence $(s_n)_{n\geq 0}$ is explicitly given by the Binet-type formula 
\begin{equation}\label{explicit formula}
    s_n = \frac{1}{2\sqrt{a^2+4b}}\left((2c_1-j_2c_0)\left(\frac{j_1}{2}\right)^n - (2c_1-j_1c_0)\left(\frac{j_2}{2}\right)^n\right), \;\;\forall n\in\mathbb N_0
\end{equation}
where $\;j_1/2 = \left(a+\sqrt{a^2+4b}\right)/2\;$ and $\;j_2/2 = \left(a-\sqrt{a^2+4b}\right)/2\;$ are roots of quadratic equation $x^2-ax-b = 0$. Both of $a$ and $b$ are positive integers, so the discriminant of $x^2-ax-b = 0$ is $a^2 + 4b > 0$. Therefore, all roots of $x^2 - ax - b =0$ are real and distinct.
\vspace{0.5cc}

Another property about $(s_n)_{n\geq 0}$ is that there exists $N\in\mathbb N$ so that $s_n \neq 0$ for all $n\geq N$. This property is useful to verify that $s_{n+1}/s_n$ will converge to a number as $\, n \,$ goes to $\infty$. If there does not exist $\: m\in\mathbb N \:$ so that $s_m = 0$, then there exists $\: N=1 \:$ so that $\: s_n \neq 0 \:$ for all $n\geq N =1$. If there exists $m\in\mathbb N$ so that $s_m = 0$, by formula \eqref{explicit formula}, we have
\begin{align*}
    (2c_1-j_2c_0)\left(\frac{j_1}{2}\right)^m = (2c_1-j_1c_0)\left(\frac{j_2}{2}\right)^m \;&\iff\; (2c_1-j_2c_0)j_1^m = (2c_1-j_1c_0)j_2^m \\
    &\iff\; (j_1^m - j_2^m)c_1 = \frac{1}{2}(j_1^{m-1}-j_2^{m-1})j_1j_2c_0 \\
    &\iff\; (j_1^m - j_2^m)c_1 = -2(j_1^{m-1}-j_2^{m-1})bc_0.
\end{align*}
If $c_0=0$ then $(j_1^m - j_2^m)c_1 = 0$. Since $j_1 \neq j_2$, then $c_1=0$. It implies $(c_0,c_1) = (0,0)$, a contradiction. Therefore $c_0$ must be non-zero and we get
\begin{equation}\label{-d/bc = xi(m)}
    \frac{-c_1}{2bc_0} = \frac{j_1^{m-1}-j_2^{m-1}}{j_1^m - j_2^m}.
\end{equation}
Let $\xi : \mathbb N \to \mathbb R$ be a function defined by $\xi(n) = \frac{j_1^{n-1}-j_2^{n-1}}{j_1^n - j_2^n}$, $\forall n\in\mathbb N$. We can observe that $\xi$ is strictly increasing and thus the equation \eqref{-d/bc = xi(m)} must have exactly one solution in $m$. Then it implies $s_n \neq 0$ for all $n\geq m+1$ (in this case, we can take $N = m+1$).
\vspace{0.3cc}
\begin{lemma}\label{s_n not 0}
    There exists $N\in\mathbb N$ so that $s_n \neq 0$ for all $n\geq N$.
\end{lemma}
\vspace{0.1cc}
\begin{lemma}
    $\: 2c_1 - j_1c_0 \:$ and $\: 2c_1-j_2c_0 \:$ can not be simultaneously equal to zero.
\end{lemma}
\begin{proof}
    Assume the contrary that $2c_1-j_1c_0 = 2c_1-j_2c_0 = 0$. By equation \eqref{explicit formula}, $s_n=0$ for all $n\in\mathbb N_0$, contradicting the Lemma \eqref{s_n not 0}. Hence, $(2c_1-j_1c_0, 2c_1-j_2c_0)\neq (0,0)$.
\end{proof}

\vspace{0.25cc}
If $2c_1-j_1c_0 = 0$, then $2c_1-j_2c_0\neq 0$ and $s_n = \frac{2c_1-j_2c_0}{2\sqrt{a^2+4b}}\left(\frac{j_1}{2}\right)^n$ for all $n\in\mathbb N_0$.

\vspace{0.3cc}
If $2c_1-j_2c_0 = 0$, then $2c_1-j_1c_0\neq 0$ and $s_n = \frac{-2c_1+j_1c_0}{2\sqrt{a^2+4b}}\left(\frac{j_2}{2}\right)^n$ for all $n\in\mathbb N_0$.

\vspace{0.3cc}
If $2c_1-j_1c_0$ and $2c_1-j_2c_0$ are non-zero: Let us consider the subsequence $(s_n)_{n\geq N}$ with $s_n\neq 0, \;\forall n\geq N$ for some $N\in\mathbb N$. By setting $n$ going to $\infty$ for $s_{n+1}/s_n$, we get
\begin{align*}
    \lim_{n\to\infty}\frac{s_{n+1}}{s_n} &= \lim_{n\to\infty}\frac{1}{2}\cdot\frac{(2c_1-j_2c_0)j_1^{n+1} - (2c_1-j_1c_0)j_2^{n+1}}{(2c_1-j_2c_0)j_1^n - (2c_1-j_1c_0)j_2^n} \\
    &= \lim_{n\to\infty}\frac{1}{2}\cdot\frac{(2c_1-j_2c_0)j_1 - (2c_1-j_1c_0)j_2\left(\frac{j_2}{j_1}\right)^n}{(2c_1-j_2c_0) - (2c_1-j_1c_0)\left(\frac{j_2}{j_1}\right)^n} \\
    &= \lim_{n\to\infty}\frac{1}{2}\cdot\frac{(2c_1-j_2c_0)j_1}{2c_1-j_2c_0} \quad\quad \left(\textrm{because $\Big|\frac{j_2}{j_1} \Big| < 1$}\right) \\
    &= \frac{j_1}{2}.
\end{align*}
Hence the sequence $\:(\frac{s_{n+1}}{s_n})_{n\geq N}\:$ converges to $\: j_1/2 = \left(a+\sqrt{a^2+4b}\right)/2$ for this case.
\vspace{0.5cc}

\begin{lemma}\label{lemma k^d s(k-1)}
    Let $d$ be a non-negative integer. So there exists the polynomials $F_{1,d}(x)$, $G_{1,d}(x)$, $H_{1,d}(x)$ with real coefficients satisfying 
    \begin{equation}\label{k^d s(k-1)}
        (\forall n\in\mathbb N) \quad 2\sum_{k=1}^{n}k^ds_{k-1} = F_{1,d}(n)s_{n+1} + G_{1,d}(n)s_n + H_{1,d}(n)
    \end{equation}
\end{lemma}

An example of ordered triple of polynomials $(F_{1,d}(x), G_{1,d}(x), H_{1,d}(x))$ that satisfies the identity \eqref{k^d s(k-1)} is stated in Theorem \eqref{matrix B and configuration of F,G,H}.
\begin{theorem}\label{matrix B and configuration of F,G,H}
    Suppose that $d\in\mathbb N_0$. Let $B_d = (b_{i,j})_{i,j=1}^{d+1}$ be a $(d+1)\times(d+1)$ square matrix with
    \begin{equation*}
       b_{i,j} = \begin{cases}
            0 & \textrm{if} \quad i>j \\
            a+b-1 & \textrm{if} \quad i=j \\
            (2^{j-i}b+a)\binom{j-1}{i-1} & \textrm{if} \quad i<j.
        \end{cases}
    \end{equation*}
    Then, we get the following statements. \\
    (i). There exists uniquely an ordered (d+1)-tuple $\left(b_0(d), b_1(d), \ldots, b_d(d)\right)\in\mathbb R^{d+1}$ satisfying the matrix equation 
    \begin{equation}\label{matrix equation B()^T = ...}
        B_d\left(b_0(d), b_1(d), \ldots, b_d(d)\right)^T = 2\left(\binom{d}{0}2^d, \binom{d}{1}2^{d-1}, \ldots, \binom{d}{d}2^0\right)^T
    \end{equation}
    (ii). Let $\left(b_0(d), b_1(d), \ldots, b_d(d)\right)$ be the unique solution of \eqref{matrix equation B()^T = ...}. \\
    Then the polynomials $F_{1,d}(x), G_{1,d}(x), H_{1,d}(x)$ defined by
    \begin{align*}
        F_{1,d}(x) = \sum_{i=0}^{d}b_i(d)x^i, \quad G_{1,d}(x) = -2(x+1)^d + \left(\sum_{i=0}^{d}(x+1)^ib_i(d)\right)b, \\
        \textrm{and} \quad H_{1,d}(x) = 2c_0 + 2^{d+1}c_1 - (ac_1+bc_0)\sum_{i=0}^{d}b_i(d) - bc_1\sum_{i=0}^{d}2^ib_i(d)
    \end{align*}
    satisfy the identity \eqref{k^d s(k-1)}.
\end{theorem}
\begin{proof}
(i). The determinant of $B_d$ is $(a+b-1)^{d+1} \neq 0$. So the matrix $B_d$ is invertible and it implies that the equation \eqref{matrix equation B()^T = ...} has exactly one solution in $\left(b_0(d), b_1(d), \ldots, b_d(d)\right)$ in $\mathbb C^{d+1}$. Furthermore,
\begin{equation*}
    \left(b_0(d), b_1(d), \ldots, b_d(d)\right)^T = 2B_d^{-1}\left(\binom{d}{0}2^d, \binom{d}{1}2^{d-1}, \ldots, \binom{d}{d}2^0\right)^T.
\end{equation*}
Since $B_d^{-1}\in\mathbb M_{d+1}(\mathbb R)$ and $\left(\binom{d}{0}2^d, \binom{d}{1}2^{d-1}, \ldots, \binom{d}{d}2^0\right)^T\in\mathbb M_{(d+1)\times 1}(\mathbb R)$, then \\ $\left(b_0(d), b_1(d), \ldots, b_d(d)\right)\in\mathbb R^{d+1}$.
\\\\
(ii). We show this part by using induction. When $n=1$, the LHS of \eqref{k^d s(k-1)} is $\: 2\sum_{k=1}^{1}k^ds_{k-1} = 2s_0 = 2c_0 \;$ and the RHS of \eqref{k^d s(k-1)} is $\; F_{1,d}(1)s_2 + G_{1,d}(1)s_1 + H_{1,d}(1) \:=\: \left(\sum_{i=0}^{d}b_i(d)\right)(ac_1+bc_0) + \left(-2^{d+1} + \left(\sum_{i=0}^{d}2^ib_i(d)\right)b\right)c_1 + 2c_0 + 2^{d+1}c_1 - (ac_1+bc_0)\sum_{i=0}^{d}b_i(d) - bc_1\sum_{i=0}^{d}2^ib_i(d) \:=\: 2c_0$. \\
In this case, the LHS and RHS of \eqref{k^d s(k-1)} are being same, so \eqref{k^d s(k-1)} satisfies for $n=1$. 

Assume that the identity \eqref{k^d s(k-1)} satisfies for $n=m$ for some $m\in\mathbb N$, so
\begin{equation}\label{induction when n=m}
    2\sum_{k=1}^{m}k^ds_{k-1} = F_{1,d}(m)s_{m+1} + G_{1,d}(m)s_m + H_{1,d}(m)
\end{equation}
We have to show that the identity \eqref{k^d s(k-1)} also holds for $n=m+1$. Before we show it, we would like to show that for all $x\in\mathbb R$,
\begin{equation}\label{ratio of F,G,H}
    F_{1,d}(x+1) : \left(F_{1,d}(x) - G_{1,d}(x+1)\right) : \left(2(x+1)^d + G_{1,d}(x)\right) = 1 : a : b
\end{equation}
,i.e., $F_{1,d}(x) - G_{1,d}(x+1) = aF_{1,d}(x+1) \;$ and $\; 2(x+1)^d + G_{1,d}(x) = bF_{1,d}(x+1)$. \\
By definitions of $F_{1,d}(x)$ and $G_{1,d}(x)$ in Theorem \eqref{matrix B and configuration of F,G,H} part (ii), we have
\begin{equation}\label{2(x+1)^d + G_1,d(x) = bF_1,d(x+1)}
    2(x+1)^d + G_{1,d}(x) = \left(\sum_{i=0}^{d}(x+1)^ib_i(d)\right)b = bF_{1,d}(x+1)
\end{equation}
Observe that the matrix equation \eqref{matrix equation B()^T = ...} is equivalent to
\begin{equation}\label{matrix equation B()^T equivalent to}
    \left(\forall k=0,1,2,\ldots,d\right) \quad \sum_{i=k}^{d}(b\cdot 2^{i-k} + a)\binom{i}{k}b_i(d) \:=\: 2\binom{d}{k}2^{d-k} + b_k(d)
\end{equation}
We also observe that the degrees of $\: F_{1,d}(x) - G_{1,d}(x+1) \:$ and $\: aF_{1,d}(x+1) \:$ are not more than $d$. For all $k=0,1,2,\ldots,d$; the coefficients of $x^k$ in polynomials $\: F_{1,d}(x) - G_{1,d}(x+1) \:$ and $\: aF_{1,d}(x+1) \:$ are respectively
\begin{equation}\label{coefficient of x^k in polynomials}
    b_k(d) + 2\binom{d}{k}2^{d-k} - \left(\sum_{i=k}^{d}2^{i-k}\binom{i}{k}b_i(d)\right)b \qquad\textrm{and}\qquad \left(\sum_{i=k}^{d}\binom{i}{k}b_i(d)\right)a
\end{equation}
Because of identity \eqref{matrix equation B()^T equivalent to}, two expressions in \eqref{coefficient of x^k in polynomials} are equivalent. It implies that
\begin{equation}\label{F_1,d(x) - G_1,d(x+1) = aF_1,d(x+1)}
    F_{1,d}(x) - G_{1,d}(x+1) = aF_{1,d}(x+1)
\end{equation}
By \eqref{2(x+1)^d + G_1,d(x) = bF_1,d(x+1)} and \eqref{F_1,d(x) - G_1,d(x+1) = aF_1,d(x+1)}, we have shown the condition \eqref{ratio of F,G,H} as desired.

We remember the recurrence relation of $(s_n)_{n\geq 0}$: $\: s_{n+1} = as_n + bs_{n-1} \:$ for all $n\in\mathbb N$. The coefficients of $s_{n+1}$, $s_n$, and $s_{n-1}$ in this recurrence relation have the ratio $\: 1:a:b \:$ which is same as the ratio in \eqref{ratio of F,G,H}. Hence, the following identity holds for all $x\in\mathbb R$ and all $n\in\mathbb N$.
\begin{equation}\label{some identity}
    F_{1,d}(x+1)s_{n+1} = \left(F_{1,d}(x)-G_{1,d}(x+1)\right)s_{n} + \left(2(x+1)^d + G_{1,d}(x)\right)s_{n-1}
\end{equation}
Setting $x:=m$, $n:=m+1$ to \eqref{some identity} yields
\begin{equation*}
    F_{1,d}(m+1)s_{m+2} = \left(F_{1,d}(m)-G_{1,d}(m+1)\right)s_{m+1} + \left(2(m+1)^d + G_{1,d}(m)\right)s_{m}
\end{equation*}
and equivalently,
\begin{equation}\label{induction, difference when n=m+1, n=m}
    2(m+1)^ds_m = F_{1,d}(m+1)s_{m+2} + \left(G_{1,d}(m+1)-F_{1,d}(m)\right)s_{m+1} - G_{1,d}(m)s_m
\end{equation}
Summing up the equations \eqref{induction when n=m} and \eqref{induction, difference when n=m+1, n=m} yields
\begin{equation*}
    2\sum_{k=1}^{m+1}k^ds_{k-1} = F_{1,d}(m+1)s_{m+2} + G_{1,d}(m+1)s_{m+1} + H_{1,d}(m+1)
\end{equation*}
and hence the identity \eqref{k^d s(k-1)} holds for $n=m+1$.

In conclusion, the identity \eqref{k^d s(k-1)} holds for all $n\in\mathbb N$.
\end{proof}
\vspace{0.5cc}

\begin{example}
    Some identities of $\:2\sum_{k=1}^{n}k^ds_{k-1}\:$ for lower degrees $d$ are:
    \begin{align*}
        &2\sum_{k=1}^{n}s_{k-1} = \frac{2s_{n+1}}{a+b-1} + \frac{(-2a+2)s_n}{a+b-1} + 2c_0 + 2c_1 - \frac{2ac_1+2bc_0+2bc_1}{a+b-1},  \\\\
        &2\sum_{k=1}^{n}ks_{k-1} = \left(\frac{2x}{a+b-1}+\frac{2a-4}{(a+b-1)^2} \right)s_{n+1} + \left(\frac{(2-2a)x}{a+b-1}+\frac{-2a^2+4a-2b-2}{(a+b-1)^2}\right)s_n \\
        & + 2c_0 + 4c_1 - \frac{(ac_1+bc_0)(4a+2b-6)+bc_1(6a+4b-8)}{(a+b-1)^2}, \\\\
        &2\sum_{k=1}^{n}k^2s_{k-1} = \left(\frac{2x^2}{a+b-1} + \frac{(4a-8)x}{(a+b-1)^2} + \frac{2a^2-6a+8b-2ab+8}{(a+b-1)^3}\right)s_{n+1} \\
        &+ \left(\frac{(2-2a)x^2}{a+b-1} + \frac{(-4a^2+8a-4b-4)x}{(a+b-1)^2} + \frac{8a^2b+2b^3+6ab^2-22ab-4b^2+18b}{(a+b-1)^3}-2\right)s_n \\
        &- \frac{(ac_1+bc_0)(8a^2+2b^2+6ab-22a-4b+18) \,+\, bc_1(18a^2+8b^2+22ab-46a-24b+32)}{(a+b-1)^3} \\
        &+ 2c_0 + 8c_1.
    \end{align*}
\end{example}
\vspace{0.5cc}

\begin{lemma}\label{zero polynomial to infinity}
    If $P(x)$ is a polynomial with real coefficients and $\;\lim\limits_{\substack{n\to\infty, \\ n\in\mathbb N}}P(n) = 0$, then $P(x)$ is zero polynomial.
\end{lemma}
\begin{proof}
    If $P(x)$ is not constant, then $\; P(x) = a_mx^m + a_{m-1}x^{m-1} + \cdots + a_1x + a_0 \:$ for some $m\in\mathbb N$ and for some $a_0, a_1, \ldots, a_m\in\mathbb R $ with $a_m\neq 0$. Observe that
    \begin{equation*}
        P(x) = a_mx^m\cdot\left(1+\sum_{k=1}^{m}\frac{a_{m-k}}{x^ka_m}\right).
    \end{equation*}
    The condition $\;\lim\limits_{n\to\infty, \, n\in\mathbb N}P(n) = 0 \;$ implies that
    \begin{equation*}
        0 = \lim\limits_{\substack{n\to\infty, \\ n\in\mathbb N}}|P(n)| = \lim\limits_{\substack{n\to\infty, \\ n\in\mathbb N}}|a_mn^m|\cdot\lim\limits_{\substack{n\to\infty, \\ n\in\mathbb N}}\Bigg|1+\sum_{k=1}^{m}\frac{a_{m-k}}{n^ka_m}\Bigg| = \lim\limits_{\substack{n\to\infty, \\ n\in\mathbb N}}|a_mn^m|.
    \end{equation*}
    It is impossible because $a_m\neq 0$ should imply $\:\lim\limits_{\substack{n\to\infty, \\ n\in\mathbb N}}|a_mn^m| = \infty$. \\
    If $P(x)$ is a constant polynomial, let $P(x) = c$, $\forall x\in\mathbb R$. Then $0 = \;\lim\limits_{\substack{n\to\infty, \\ n\in\mathbb N}}P(n) = \;\lim\limits_{\substack{n\to\infty, \\ n\in\mathbb N}} c = c$, so $P(x)$ is zero polynomial. \\
    In conclusion, $P(x)$ must be zero polynomial.
\end{proof}

\vspace{1cc}
\begin{lemma}\label{fraction x^t per k^x}
    Let $t$ be a positive integer and $k>1$ be a real number. Then
    \begin{equation*}
         \:\lim\limits_{\substack{n\to\infty \\ n\in\mathbb N}}\frac{n^t}{k^n} = 0.
    \end{equation*}
\end{lemma}
\begin{proof}
    For all positive integers $n > \frac{(t+2)!}{(\ln k)^{t+2}}$, we have
    \begin{equation*}
        0 < \frac{n^{t+1}}{k^n} = \frac{n^{t+1}}{\sum_{a=0}^{\infty}\frac{(n\ln k)^a}{a!}} < \frac{n^{t+1}}{\frac{(n\ln k)^{t+2}}{(t+2)!}} = \frac{(t+2)!}{n(\ln k)^{t+2}} < 1 \;\implies\; 0 < \frac{n^t}{k^n} < \frac{1}{n}.
    \end{equation*}
    Since $\; \lim\limits_{\substack{n\to\infty \\ n\in\mathbb N}} 0 = \lim\limits_{\substack{n\to\infty \\ n\in\mathbb N}}\frac{1}{n} = 0$, by Squeeze theorem we obtain $\; \lim\limits_{\substack{n\to\infty \\ n\in\mathbb N}}\frac{n^t}{k^n} = 0$. 
\end{proof}

\vspace{1cc}
\begin{theorem}\label{polinom gamma1, gamma2, gamma3}
    Let $\gamma_1(x)$, $\gamma_2(x)$, $\gamma_3(x)$ be polynomials with real coefficients. Let $2c_1-j_1c_0$ and $2c_1-j_2c_0$ be non-zero. If $\: \gamma_1(n)s_{n+1} + \gamma_2(n)s_n + \gamma_3(n) = 0 \:$ holds for all $n\in\mathbb N$, then $\gamma_1(x)$, $\gamma_2(x)$, $\gamma_3(x)$ are zero polynomials.
\end{theorem}
\begin{proof}
    We know from Lemma \eqref{s_n not 0} that there exists a number $N\in\mathbb N$ so that $s_n \neq 0$ for all $n\geq N$. Dividing both sides of identity $\: \gamma_1(n)s_{n+1} + \gamma_2(n)s_n + \gamma_3(n) = 0 \:$ by $s_n$ over integers $n\geq N$ implies
    \begin{equation}\label{gamma1,2,3 when n>=N}
        \left(\forall n\in\mathbb N, \: n\geq N\right) \;\; \gamma_1(n)\frac{s_{n+1}}{s_n} + \gamma_2(n) + \frac{\gamma_3(n)}{s_n} = 0
    \end{equation}
    Observe that $j_1/2 = \frac{a+\sqrt{a^2+4b}}{2} > a \geq 1$. By using Lemma \eqref{fraction x^t per k^x} and the explicit formula \eqref{explicit formula} of $s_n$, we should have
    \begin{equation*}
        \lim\limits_{\substack{n\to\infty \\ n\in\mathbb N}}\frac{n^{\deg(\gamma_3)+1}}{s_n} = \frac{2\sqrt{a^2+4b}}{2c_1-j_2c_0}\cdot\lim\limits_{\substack{n\to\infty \\ n\in\mathbb N}}\frac{n^{\deg(\gamma_3)+1}}{\left(\frac{j_1}{2}\right)^n} = 0
    \end{equation*}
    where $\deg(...)$ is degree of a polynomial. Therefore,
    \begin{equation*}
        \lim\limits_{n\to\infty}\frac{\gamma_3(n)}{s_n} = \lim\limits_{n\to\infty}\frac{\gamma_3(n)}{n^{\deg(\gamma_3)+1}}\cdot\lim\limits_{n\to\infty}\frac{n^{\deg(\gamma_3)+1}}{s_n} = 0\cdot 0 = 0.
    \end{equation*}
    By setting $n$ going to $\infty$ in \eqref{gamma1,2,3 when n>=N}, it yields $\;\lim\limits_{\substack{n\to\infty \\ n\in\mathbb N}}\left(\gamma_1(n)\frac{j_1}{2} + \gamma_2(n)\right) = 0. \;$ \\
    Since $\;\gamma_1(x)\frac{j_1}{2} + \gamma_2(x)\in\mathbb R[x]$, by Lemma \eqref{zero polynomial to infinity}, then $\;\gamma_1(x)\frac{j_1}{2} + \gamma_2(x)\;$ is a zero polynomial. Substituting $\:\gamma_2(x) = -\frac{j_1}{2}\gamma_1(x), \forall x\in\mathbb R\;$ to identity $\: \gamma_1(n)s_{n+1} + \gamma_2(n)s_n + \gamma_3(n) = 0, \forall n\in\mathbb N\:$ yields
    \begin{align*}
       \left(\forall n\in\mathbb N\right)\;\; 0 = \gamma_1(n)s_{n+1} + \gamma_2(n)s_n + \gamma_3(n) 
       &= \gamma_1(n)\cdot\left(s_{n+1}-\frac{j_1}{2}s_n\right) + \gamma_3(n) \\
       &= \frac{2c_1-j_1c_0}{2}\left(\frac{j_2}{2}\right)^n\gamma_1(n) + \gamma_3(n) \quad
    \end{align*}
    \begin{equation}\label{gamma1 and gamma3 equation}
        \implies\;\;\;\; (\forall n\in\mathbb N)\;\;\;\frac{2c_1-j_1c_0}{2}\left(\frac{j_2}{2}\right)^n\gamma_1(n) + \gamma_3(n) = 0 \qquad\quad
    \end{equation}
    If one of $\gamma_1(x)$ or $\gamma_3(x)$ is a zero polynomial, by equation \eqref{gamma1 and gamma3 equation} it implies that both of $\gamma_1(x)$ and $\gamma_3(x)$ are zero polynomials. Then $\gamma_2$ is a zero polynomial too, and the theorem is done. \\
    If both of $\gamma_1$ and $\gamma_3$ are not zero polynomials, then there exists $M_1\in\mathbb N$ such that $\gamma_1(x)\neq 0$ and $\gamma_3(x)\neq 0$ for all $x\geq M_1$. From the identity \eqref{gamma1 and gamma3 equation}, we have
    \begin{equation}\label{frac polynom gamma3 per gamma1}
        \left(\forall n\in\mathbb N, n\geq M_1\right)\quad \frac{\gamma_3(n)}{\gamma_1(n)} = \frac{-(2c_1-j_1c_0)}{2}\left(\frac{j_2}{2}\right)^n \qquad
    \end{equation}
    then
    \begin{equation}\label{absolute frac gamma3 per gamma1}
        \left(\forall n\in\mathbb N, n\geq M_1\right)\quad \bigg|\frac{\gamma_3(n)}{\gamma_1(n)}\bigg| = \frac{1}{2}|2c_1-j_1c_0|\bigg|\frac{j_2}{2}\bigg|^n \qquad
    \end{equation}
    Let us observe \eqref{frac polynom gamma3 per gamma1} and \eqref{absolute frac gamma3 per gamma1} into 3 possibilities: $b-a>1$, $b-a=1$, and $b-a<1$. \\\\
    \underline{Case $b-a>1$}: \\
    Observe that $|j_2/2| = \frac{\sqrt{a^2+4b}-a}{2}>1$ and
    \begin{equation*}
        \lim\limits_{\substack{n\to\infty \\ n\in\mathbb N}}\frac{|\gamma_3(n)/\gamma_1(n)|}{n^{|\deg(\gamma_3)-\deg(\gamma_1)|+1}} = 0.
    \end{equation*}
    By Lemma \eqref{fraction x^t per k^x}, we also have
    \begin{equation*}
        \lim\limits_{\substack{n\to\infty \\ n\in\mathbb N}}\frac{n^{|\deg(\gamma_3)-\deg(\gamma_1)|+1}}{|j_2/2|^n} = 0.
    \end{equation*}
    Therefore, from \eqref{absolute frac gamma3 per gamma1},
    \begin{align*}
        |2c_1-j_1c_0| &= 2\lim\limits_{\substack{n\to\infty \\ n\in\mathbb N}}\frac{|\gamma_3(n)/\gamma_1(n)|}{|j_2/2|^n} \\
        &= \lim\limits_{\substack{n\to\infty \\ n\in\mathbb N}}\frac{|\gamma_3(n)/\gamma_1(n)|}{n^{|\deg(\gamma_3)-\deg(\gamma_1)|+1}}\cdot\lim\limits_{\substack{n\to\infty \\ n\in\mathbb N}}\frac{n^{|\deg(\gamma_3)-\deg(\gamma_1)|+1}}{|j_2/2|^n} \\
        &= 0\cdot 0 = 0
    \end{align*}
    then $\: 2c_1-j_1c_0 = 0$, a contradiction. \\\\
    \underline{Case $b-a = 1$}: \\
    We now have $j_2/2 = \frac{a-\sqrt{a^2+4b}}{2} = -1$, so the identities \eqref{frac polynom gamma3 per gamma1} and \eqref{absolute frac gamma3 per gamma1} respectively become
    \begin{equation}\label{first main eq when b-a=1}
         \left(\forall n\in\mathbb N, n\geq M_1\right)\quad \frac{\gamma_3(n)}{\gamma_1(n)} = \frac{-(2c_1-j_1c_0)}{2}\left(-1\right)^n \qquad
    \end{equation}
    and
    \begin{equation}\label{second main eq when b-a=1}
        \left(\forall n\in\mathbb N, n\geq M_1\right)\quad \bigg|\frac{\gamma_3(n)}{\gamma_1(n)}\bigg| = \frac{1}{2}|2c_1-j_1c_0| \qquad
    \end{equation}
    From \eqref{second main eq when b-a=1}, it indicates that either $\:\gamma_3(n) = \frac{1}{2}(2c_1-j_1c_0)\gamma_1(n)\:$ for infinitely many integers $n$ or $\:\gamma_3(n) = -\frac{1}{2}(2c_1-j_1c_0)\gamma_1(n)\:$ for infinitely many integers $n$. \\
    If $\:\gamma_3(n) = \frac{1}{2}(2c_1-j_1c_0)\gamma_1(n)\:$ for infinitely many integers $n$, it implies $\: \gamma_3(x) = \frac{1}{2}(2c_1-j_1c_0)\gamma_1(x)\:$ for all $x\in\mathbb R$. Therefore, we can rewrite \eqref{first main eq when b-a=1} by
    \begin{equation}\label{third eq when b-a=1}
        \left(\forall n\in\mathbb N, n\geq M_1\right)\quad \frac{2c_1-j_1c_0}{2} = \frac{-(2c_1-j_1c_0)}{2}\left(-1\right)^n \qquad
    \end{equation}
    Setting to \eqref{third eq when b-a=1} when $n$ is even yields that $\: 2c_1-j_1c_0 = 0$, a contradiction. \\
    If $\:\gamma_3(n) = -\frac{1}{2}(2c_1-j_1c_0)\gamma_1(n)\:$ for infinitely many integers $n$, then $\:\gamma_3(x) = -\frac{1}{2}(2c_1-j_1c_0)\gamma_1(x)\:$ for all $x\in\mathbb R$. So we can rewrite  \eqref{first main eq when b-a=1} by
    \begin{equation}\label{fourth eq when b-a=1}
        \left(\forall n\in\mathbb N, n\geq M_1\right)\quad -\frac{1}{2}(2c_1-j_1c_0) = \frac{-(2c_1-j_1c_0)}{2}\left(-1\right)^n \qquad
    \end{equation}
    Setting to \eqref{fourth eq when b-a=1} when $n$ is odd yields that $2c_1-j_1c_0 = 0$, a contradiction. \\\\
    \underline{Case $b-a < 1$}: \\
    We have that $0 < |j_2/2| = \frac{\sqrt{a^2+4b}-a}{2} < 1$, therefore $|j_2/2| > 1$. By Lemma \eqref{fraction x^t per k^x}, 
    \begin{equation*}
        \lim\limits_{\substack{n\to\infty \\ n\in\mathbb N}}\frac{n^{|\deg(\gamma_3)-\deg(\gamma_1)|+1}}{|2/j_2|^n} = 0.
    \end{equation*}
    Besides that,
    \begin{equation*}
        \lim\limits_{\substack{n\to\infty \\ n\in\mathbb N}}\frac{|\gamma_1(n)/\gamma_3(n)|}{n^{|\deg(\gamma_3)-\deg(\gamma_1)|+1}} = 0.
    \end{equation*}
    Hence, from \eqref{absolute frac gamma3 per gamma1},
    \begin{equation*}
        \bigg|\frac{2}{2c_1-j_1c_0}\bigg| = \lim\limits_{\substack{n\to\infty \\ n\in\mathbb N}}\frac{|\gamma_1(n)/\gamma_3(n)|}{|2/j_2|^n} 
        = \lim\limits_{\substack{n\to\infty \\ n\in\mathbb N}}\frac{|\gamma_1(n)/\gamma_3(n)|}{n^{|\deg(\gamma_3)-\deg(\gamma_1)|+1}}\cdot\lim\limits_{\substack{n\to\infty \\ n\in\mathbb N}}\frac{n^{|\deg(\gamma_3)-\deg(\gamma_1)|+1}}{|2/j_2|^n} 
        = 0\cdot 0 = 0
    \end{equation*}
    which is impossible since $\big|\frac{2}{2c_1-j_1c_0}\big|$ itself is non-zero. \\\\
    In conclusion, $\gamma_1(x)$, $\gamma_2(x)$, $\gamma_3(x)$ must be zero polynomials.
\end{proof}
\vspace{1cc}
\begin{theorem}\label{other two conditions for gamma1, gamma2, gamma3}
    Let $\gamma_1(x)$, $\gamma_2(x)$, $\gamma_3(x)$ be polynomials with real coefficients, $0\in\{2c_1-j_1c_0, 2c_1-j_2c_0\}$, and $\gamma_1(n)s_{n+1}+\gamma_2(n)s_n + \gamma_3(n) = 0$ for all $n\in\mathbb N$. \\
    $(i)$. If $\: 2c_1-j_1c_0 = 0$, then $\: \gamma_2 \equiv \frac{-j_1}{2}\cdot\gamma_1 \:$ and $\: \gamma_3 \equiv 0$. \\
    $(ii)$. If $\: 2c_1-j_2c_0 = 0$, then $\: \gamma_2 \equiv \frac{-j_2}{2}\cdot\gamma_1 \:$ and $\: \gamma_3 \equiv 0$.
\end{theorem}
\begin{proof}
    (i). If $\: 2c_1-j_1c_0 = 0$, then $s_n = (j_1/2)^nK_1$ for all $n\in\mathbb N_0$, where $K_1 = \frac{2c_1-j_2c_0}{2\sqrt{a^2+4b}}$ is a non-zero real constant. The identity $\:\gamma_1(n)s_{n+1}+\gamma_2(n)s_n + \gamma_3(n) = 0, \:\forall n\in\mathbb N \:$ becomes
    \begin{equation}\label{equation case 2c1-j1c0 = 0}
        \left(\gamma_1(n)\frac{j_1}{2} + \gamma_2(n)\right)\left(\frac{j_1}{2}\right)^nK_1 + \gamma_3(n) = 0, \;\;\forall n\in\mathbb N 
    \end{equation}

    \begin{equation*}    
        \iff\;\;  \gamma_1(n)\frac{j_1}{2} + \gamma_2(n) = \frac{-1}{K_1}\cdot\frac{\gamma_3(n)}{\left(j_1/2\right)^n}, \;\;\forall n\in\mathbb N.
    \end{equation*}
    Observe that
    \begin{equation*}
        \lim\limits_{\substack{n\to\infty \\ n\in\mathbb N}}\frac{\gamma_3(n)}{n^{\deg(\gamma_3)+1}} = 0
    \end{equation*}
    and $j_1/2 > 1$. By Theorem \eqref{fraction x^t per k^x}, we have
    \begin{equation*}
        \lim\limits_{\substack{n\to\infty \\ n\in\mathbb N}}\frac{n^{\deg(\gamma_3)+1}}{(j_1/2)^n} = 0.
    \end{equation*}
    Consequently,
    \begin{align*}
        \lim\limits_{\substack{n\to\infty \\ n\in\mathbb N}}\left(\gamma_1(n)\frac{j_1}{2} + \gamma_2(n)\right) = \lim\limits_{\substack{n\to\infty \\ n\in\mathbb N}}\frac{-1}{K_1}\cdot\frac{\gamma_3(n)}{\left(j_1/2\right)^n} &= \frac{-1}{K_1}\cdot\lim\limits_{\substack{n\to\infty \\ n\in\mathbb N}}\frac{\gamma_3(n)}{n^{\deg(\gamma_3)+1}}\cdot\lim\limits_{\substack{n\to\infty \\ n\in\mathbb N}}\frac{n^{\deg(\gamma_3)+1}}{(j_1/2)^n} \\
        &= \frac{-1}{K_1}\cdot 0\cdot 0 = 0.
    \end{align*}
    Since $\gamma_1(n)\frac{j_1}{2} + \gamma_2(n)\in\mathbb R[x]$, by Theorem \eqref{zero polynomial to infinity}, it leads to the result that $\frac{j_1}{2}\cdot\gamma_1 + \gamma_2$ is zero polynomial, then $\gamma_2 \equiv \frac{-j_1}{2}\gamma_1$. Also, equation \eqref{equation case 2c1-j1c0 = 0} implies $\gamma_3 \equiv 0$. \\\\
    (ii). If $\: 2c_1-j_2c_0 = 0$, we have $s_n = (j_2/2)^nK_2$ for all $n\in\mathbb N$, where $K_2 = \frac{-2c_1+j_1c_0}{2\sqrt{a^2+4b}}$ is a non-zero real constant. Therefore, we can rewrite the identity $\gamma_1(n)s_{n+1}+\gamma_2(n)s_n + \gamma_3(n) = 0, \:\forall n\in\mathbb N$ by
    \begin{equation}\label{eq1 case 2c1-j2c0 = 0}
        \left(\gamma_1(n)\frac{j_2}{2} + \gamma_2(n)\right)\left(\frac{j_2}{2}\right)^nK_2 + \gamma_3(n) = 0, \;\;\forall n\in\mathbb N
    \end{equation}
    Let us denote the polynomial $\gamma_1(x)\frac{j_2}{2} + \gamma_2(x)$ by $\gamma_4(x)$. The identity \eqref{eq1 case 2c1-j2c0 = 0} is equivalent to
    \begin{equation}\label{eq2 case 2c1-j2c0 = 0}
        \gamma_4(n)\left(\frac{j_2}{2}\right)^nK_2 + \gamma_3(n) = 0, \;\;\forall n\in\mathbb N
    \end{equation}
    Note that if either $\gamma_3$ or $\gamma_4$ is zero polynomial, then both of them becomes zero polynomials. Hence $\: \gamma_2 \equiv \frac{-j_2}{2}\cdot\gamma_1 \:$ and $\: \gamma_3 \equiv 0$, done. \\\\
    If $\gamma_3$ and $\gamma_4$ are not zero polynomials, then there exists $M_2\in\mathbb N$ in such a way that $\gamma_3(x)\neq 0$, $\gamma_4(x)\neq 0$ for all real $x\geq M_2$. So the identity \eqref{eq2 case 2c1-j2c0 = 0} implies the following two identities:
    \begin{equation}\label{eq3 case 2c1-j2c0 = 0}
        \left(\forall n\in\mathbb N, n\geq M_2\right)\quad \frac{\gamma_3(n)}{\gamma_4(n)} = \frac{-K_2}{2}\left(\frac{j_2}{2}\right)^n
    \end{equation}
    \begin{equation}\label{eq4 case 2c1-j2c0 = 0}
        \left(\forall n\in\mathbb N, n\geq M_2\right)\quad \bigg|\frac{\gamma_3(n)}{\gamma_4(n)}\bigg| = \frac{|K_2|}{2}\bigg|\frac{j_2}{2}\bigg|^n
    \end{equation}
    We intend to investigate \eqref{eq3 case 2c1-j2c0 = 0} and \eqref{eq4 case 2c1-j2c0 = 0} in $3$ cases: $b>a+1$, $\: b=a+1$, $\: b<a+1$. \\
    \underline{Case $b>a+1$}: We have $|j_2/2| > 1$ and
    \begin{equation*}
        \lim\limits_{\substack{n\to\infty \\ n\in\mathbb N}}\frac{|\gamma_3(n)/\gamma_4(n)|}{n^{|\deg(\gamma_3)-\deg(\gamma_4)|+1}} = 0.
    \end{equation*}
    By Theorem \eqref{fraction x^t per k^x}, we get
    \begin{equation*}
        \lim\limits_{\substack{n\to\infty \\ n\in\mathbb N}}\frac{n^{|\deg(\gamma_3)-\deg(\gamma_4)|+1}}{|j_2/2|^n} = 0.
    \end{equation*}
    Hence by \eqref{eq4 case 2c1-j2c0 = 0}, it leads to the result
    \begin{equation*}
        \frac{|K_2|}{2} = \lim\limits_{\substack{n\to\infty \\ n\in\mathbb N}}\frac{|\gamma_3(n)/\gamma_4(n)|}{|j_2/2|^n} = \lim\limits_{\substack{n\to\infty \\ n\in\mathbb N}}\frac{|\gamma_3(n)/\gamma_4(n)|}{n^{|\deg(\gamma_3)-\deg(\gamma_4)|+1}} \cdot  \lim\limits_{\substack{n\to\infty \\ n\in\mathbb N}}\frac{n^{|\deg(\gamma_3)-\deg(\gamma_4)|+1}}{|j_2/2|^n} = 0\cdot 0 = 0,
    \end{equation*}
    implying that $K_2 = 0$, it is a contradiction.
    \\\\
    \underline{Case $b=a+1$}: 
    In this case, $j_2/2 = 1$. The identities \eqref{eq3 case 2c1-j2c0 = 0} and \eqref{eq4 case 2c1-j2c0 = 0} can be rewritten as the following two:
    \begin{equation}\label{eq5 case 2c1-j2c0 = 0}
        \left(\forall n\in\mathbb N, n\geq M_2\right)\quad \frac{\gamma_3(n)}{\gamma_4(n)} = \frac{-K_2}{2}\left(-1\right)^n
    \end{equation}
    \begin{equation}\label{eq6 case 2c1-j2c0 = 0}
        \left(\forall n\in\mathbb N, n\geq M_2\right)\quad \bigg|\frac{\gamma_3(n)}{\gamma_4(n)}\bigg| = \frac{|K_2|}{2}
    \end{equation}
    The identity \eqref{eq6 case 2c1-j2c0 = 0} gives us the fact that either $\gamma_3(n) = \frac{1}{2}K_2\gamma_4(n)$ for infinitely many integers $n$ or $\gamma_3(n) = -\frac{1}{2}K_2\gamma_4(n)$ for infinitely many integers $n$.

    If $\:\gamma_3(n) = \frac{1}{2}K_2\gamma_4(n)\:$ for infinitely many integers $n$, it implies that $\:\gamma_3(x) = \frac{1}{2}K_2\gamma_4(x)\:$ for all real $x$. Substituting it to \eqref{eq5 case 2c1-j2c0 = 0} and setting when $n$ is even yields that $K_2 = 0$, a contradiction.

    If $\:\gamma_3(n) = -\frac{1}{2}K_2\gamma_4(n)\:$ for infinitely many integers $n$, then $\:\gamma_3(x) = -\frac{1}{2}K_2\gamma_4(x)\:$ for all real $x$. Substituting it to \eqref{eq5 case 2c1-j2c0 = 0} and setting when $n$ is odd yields that $K_2 = 0$, a contradiction.
    \\\\
    \underline{Case $b<a+1$}: In this case, we have the fact $0 < |j_2/2| < 1$, so $|2/j_2| > 1$. Notice that 
    \begin{equation*}
        \lim\limits_{\substack{n\to\infty \\ n\in\mathbb N}}\frac{|\gamma_4(n)/\gamma_3(n)|}{n^{|\deg(\gamma_4)-\deg(\gamma_3)|+1}} = 0
    \end{equation*}
    and by Theorem \eqref{fraction x^t per k^x}, we also have
    \begin{equation*}
        \lim\limits_{\substack{n\to\infty \\ n\in\mathbb N}}\frac{n^{|\deg(\gamma_4)-\deg(\gamma_3)|+1}}{|2/j_2|^n} = 0.
    \end{equation*}
    By \eqref{eq4 case 2c1-j2c0 = 0}, we obtain
    \begin{equation*}
        \frac{2}{|K_2|} = \lim\limits_{\substack{n\to\infty \\ n\in\mathbb N}}\frac{|\gamma_4(n)/\gamma_3(n)|}{|2/j_2|^n} = \lim\limits_{\substack{n\to\infty \\ n\in\mathbb N}}\frac{|\gamma_4(n)/\gamma_3(n)|}{n^{|\deg(\gamma_4)-\deg(\gamma_3)|+1}}\cdot\lim\limits_{\substack{n\to\infty \\ n\in\mathbb N}}\frac{n^{|\deg(\gamma_4)-\deg(\gamma_3)|+1}}{|2/j_2|^n} = 0\cdot 0 = 0,
    \end{equation*}
    but it is impossible. \\\\
    In conclusion, the result must be $\: \gamma_2 \equiv \frac{-j_2}{2}\cdot\gamma_1 \:$ and $\: \gamma_3 \equiv 0$. 
\end{proof}

\vspace{1.5cc}

\section{Main Theorems}\label{section 3}
\noindent
This section presents three main theorems of this paper. Theorem \eqref{main theorem 1} is a development of Lemma \eqref{lemma k^d s(k-1)} and Theorem \eqref{matrix B and configuration of F,G,H}, because Lemma \eqref{lemma k^d s(k-1)} and Theorem \eqref{matrix B and configuration of F,G,H} only represent the summation $2\sum_{k=1}^{n}k^ds_{k-1}$ for arbitrary $d\in\mathbb N_0$, while Theorem \eqref{main theorem 1} represents the summation $2\sum_{k=1}^{n}P(k)s_{k-1}$ for all general polynomials $P(x)\in\mathbb R[x]$. For an arbitrary $P(x)\in\mathbb R[x]$, we can find the existence of polynomials $F_1(x), G_1(x), H_1(x)\in\mathbb R[x]$ satisfying the identity $\; 2\sum_{k=1}^{n}P(k)s_{k-1} = F_1(n)s_{n+1} + G_1(n)s_n + H_1(n), \;\forall n\in\mathbb N$. Not only that, we also determine how many triples $(F_1(x), G_1(x), H_1(x))$ satisfying that identity. Here, we get the result that if both of $2c_1-j_1c_0$ and $2c_1-j_2c_0$ are non-zero then the triple $(F_1(x), G_1(x), H_1(x))$ is unique. If one of $2c_1 - j_1c_0$ or $2c_1-j_2c_0$ equals zero, so there are infinitely many triples $(F_1(x), G_1(x), H_1(x))$.

\vspace{0.5cc}
\begin{theorem}\label{main theorem 1}
    Let $P(x)$ be a polynomial with real coefficients. Then there exists a triple of polynomials $(F_1(x), G_1(x), H_1(x))$ with real coefficients so that 
    \begin{equation}\label{main identity in main theorem}
        2\sum_{k=1}^{n}P(k)s_{k-1} = F_1(n)s_{n+1} + G_1(n)s_n + H_1(n), \;\forall n\in\mathbb N.
    \end{equation}
\end{theorem}
\begin{proof}
    Consider a polynomial $P(x)\in\mathbb R[x]$. \\
    We can state that $P(x) = a_mx^m + a_{m-1}x^{m-1} + \cdots + a_1x + a_0$ for some $m\in\mathbb N_0$ and $a_0, a_1, \cdots, a_m\in\mathbb R$. An example of triple $(F_1(x), G_1(x), H_1(x))$ satisfying \eqref{main identity in main theorem} is
    \begin{equation*}
        F_1(x) = \sum_{d=0}^{m}a_dF_{1,d}(x); \quad G_1(x) = \sum_{d=0}^{m}a_dG_{1,d}(x); \quad H_1(x) = \sum_{d=0}^{m}a_dH_{1,d}(x) 
    \end{equation*}
    where $F_{1,d}(x)$, $G_{1,d}(x)$, $H_{1,d}(x)$ are polynomials as defined in Theorem \eqref{matrix B and configuration of F,G,H} part (ii). \\
    Hence, the result follows.
\end{proof}

\vspace{0.5cc}
\begin{theorem}\label{main theorem 2}
    If both of $2c_1-j_1c_0$ and $2c_1-j_2c_0$ are non-zero, then for every polynomials $P(x)\in\mathbb R[x]$, the triple $(F_1(x), G_1(x), H_1(x))\in\mathbb R[x]^3$ satisfying \eqref{main identity in main theorem} is unique.
\end{theorem}
\begin{proof}
    Suppose that $(F_2(x), G_2(x), H_2(x))$ and $(F_3(x), G_3(x), H_3(x))$ are the ordered solution of $(F_1(x), G_1(x), H_1(x))$ which satisfy \eqref{main identity in main theorem}. Then we have
    \begin{align*}
        2\sum_{k=1}^{n}P(k)s_{k-1} = F_2(n)s_{n+1} + G_2(n)s_n + H_2(n) = F_3(n)s_{n+1} + G_3(n)s_n + H_3(n), \;\;\forall n\in\mathbb N
    \end{align*}
    then
    \begin{equation*}
        (F_2-F_3)(n)s_{n+1} + (G_2-G_3)(n)s_n + (H_2-H_3)(n) = 0, \;\; \forall n\in\mathbb N.
    \end{equation*}
    Since $(F_2-F_3)(x)$, $(G_2-G_3)(x)$, $(H_2-H_3)(x)$ are polynomials with real coefficients, by Theorem \eqref{polinom gamma1, gamma2, gamma3} we get that $(F_2-F_3)(x)$, $(G_2-G_3)(x)$, $(H_2-H_3)(x)$ are identically zero. Therefore $F_2(x)$, $G_2(x)$, $H_2(x)$ are identically equal to $F_3(x)$, $G_3(x)$, $H_3(x)$ respectively. \\
    Hence the triple of polynomials $(F_1(x), G_1(x), H_1(x))$ which satisfies \eqref{main identity in main theorem} is unique.
\end{proof}

\vspace{0.5cc}
\begin{theorem}\label{main theorem 3}
    If either $2c_1-j_1c_0$ or $2c_1-j_2c_0$ is equal to $0$, then for every polynomials $P(x)\in\mathbb R[x]$, the triple $(F_1(x), G_1(x), H_1(x))\in\mathbb R[x]^3$ satisfying \eqref{main identity in main theorem} is infinitely many. 
\end{theorem}
\begin{proof}
    Define that $\Omega_1$ and $\Omega_2$ are the sets of all triples $(F_1(x), G_1(x), H_1(x))\in\mathbb R[x]^3$ satisfying \eqref{main identity in main theorem} with the constraints $2c_0-j_1c_0 = 0$ and $2c_1-j_2c_0 = 0$ respectively. Define that $P(x) = a_mx^m + a_{m-1}x^{m-1} + \cdots + a_1x + a_0$ for some $m\in\mathbb N_0$ and $a_0, a_1, \cdots, a_m\in\mathbb R$.

    If $\: 2c_1-j_1c_0 = 0$, by Theorem \eqref{other two conditions for gamma1, gamma2, gamma3}, we have the conditions that for all $(F_2, G_2, H_2)$ and $(F_3, G_3, H_3)$ in $\Omega_1$,
    \begin{equation*}
        \frac{-j_1}{2}\cdot (F_2 - F_3) \equiv G_2 - G_3 \quad\textrm{and}\quad H_2 - H_3 \equiv 0
    \end{equation*}
    \begin{equation*}
        \iff\quad \frac{j_1}{2}\cdot F_2 + G_2 \equiv \frac{j_1}{2}\cdot F_3 + G_3(x) \quad\textrm{and}\quad H_2 \equiv H_3.
    \end{equation*}
    We know that $\left(\sum_{d=0}^{m}a_dF_{1,d}(x), \: \sum_{d=0}^{m}a_dG_{1,d}(x), \: \sum_{d=0}^{m}a_dH_{1,d}(x)\right)$ is an element in $\Omega_1$, so it leads to the fact that for all $(F, G, H)$ in $\Omega_1$,
    \begin{equation*}
        \frac{j_1}{2}\cdot F + G \equiv \frac{j_1}{2}\sum_{d=0}^{m}a_dF_{1,d} + \sum_{d=0}^{m}a_dG_{1,d} \quad\textrm{and}\quad H \equiv \sum_{d=0}^{m}a_dH_{1,d}.
    \end{equation*}
    Therefore
    \begin{equation*}
        \Omega_1 = \{F(x), \frac{-j_1}{2}\cdot F(x) + \frac{j_1}{2}\sum_{d=0}^{m}a_dF_{1,d}(x) + \sum_{d=0}^{m}a_dG_{1,d}(x), \sum_{d=0}^{m}a_dH_{1,d}(x) \:\big | \: F(x)\in\mathbb R[x]\}
    \end{equation*}
    and it is obvious that $|\Omega_1| = \infty$. 

    It is also similar when $2c_1-j_2c_0 = 0$. We will directly get
    \begin{equation*}
        \Omega_2 = \{F(x), \frac{-j_2}{2}\cdot F(x) + \frac{j_2}{2}\sum_{d=0}^{m}a_dF_{1,d}(x) + \sum_{d=0}^{m}a_dG_{1,d}(x), \sum_{d=0}^{m}a_dH_{1,d}(x) \:\big | \: F(x)\in\mathbb R[x]\}.
    \end{equation*}
    and then $|\Omega_2| = \infty$. \\
    Hence, $|\Omega_1| = |\Omega_2| = \infty \:$ answers the Theorem.
\end{proof}

\bigskip
\hrule
\bigskip

\noindent 2020 {\it Mathematics Subject Classification}:
Primary 11B39; Secondary 11B83.

\noindent \emph{Keywords}: generalized Fibonacci numbers, polynomials, finite series.

\bigskip
\hrule
\bigskip

\noindent (Concerned with sequence
\seqnum{A000045}, \seqnum{A000129}, \seqnum{A000032}, \seqnum{A001045}, and many others.)

\end{document}